\RequirePackage{fix-cm}
\documentclass[smallextended,draft]{svjour3}       
\smartqed  
\usepackage{graphicx}
\usepackage[T1]{fontenc}
\usepackage[utf8]{inputenc}

\usepackage{amsmath,amssymb,amsfonts}
\usepackage{yfonts}
\usepackage[dvipsnames]{xcolor}
\usepackage{cite}
\usepackage{enumitem}

 \usepackage{newtxtext}      

\newtheorem{assumption}{Assumption}

\newcommand{\RR}{\mathbb R}
\newcommand{\TT}{\mathbb T}

\newcommand{\cA}{\mathcal{A}}
\newcommand{\cB}{\mathcal{B}}
\newcommand{\cC}{\mathcal{C}}
\newcommand{\cL}{\mathcal{L}}

\DeclareMathOperator{\sat}{sat}

\DeclareMathOperator{\rank}{rank}
\DeclareMathOperator{\diag}{diag}
\DeclareMathOperator{\Ran}{Ran}
\DeclareMathOperator{\Ker}{Ker}
\def\d{\mathrm{d}}
\def\dt{\mathrm{d} t}
\def\dx{\mathrm{d} x}

\newcommand{\overbar}[1]{\mkern 1.5mu\overline{\mkern-1.5mu#1\mkern-1.5mu}\mkern 1.5mu}

 \journalname{Mathematics of Control, Signals, and Systems}
\begin{document}

\title{Forwarding techniques for the global stabilization of dissipative infinite-dimensional systems coupled with an ODE\thanks{This research was partially supported by the French Grant ANR ODISSE (ANR-19-CE48-0004-01) and was also conducted in the framework of the regional programme "Atlanstic 2020, Research, Education and Innovation in Pays de la Loire“, supported by the French Region Pays de la Loire and the European Regional Development Fund.}
}

\titlerunning{Forwarding techniques for the global stabilization of dissipative PDEs coupled with an ODE}        

\author{ Swann Marx \and Lucas Brivadis \and  Daniele Astolfi        
}


\institute{
S. Marx \at
        LS2N, \'Ecole Centrale de Nantes $\&$ CNRS UMR 6004, F-44000 Nantes, France.    \\
              \email{swann.marx@ls2n.fr}           
\and 
L. Brivadis and D. Astolfi \at
              Univ Lyon, Universit\'e Claude Bernard Lyon 1, CNRS, LAGEPP UMR 5007, 43 boulevard du 11 novembre 1918, F-69100, Villeurbanne, France \\
              \email{lucas.brivadis@univ-lyon1.fr; daniele.astolfi@univ-lyon1.fr}
}

\date{Received: \today / Accepted: date}

\maketitle

\begin{abstract}
This paper deals  with the stabilization  of a coupled system composed by an infinite-dimensional system and an ODE. Moreover, the control, which appears in the dynamics of the ODE, is subject to a general class of nonlinearities. Such a situation may arise, for instance, 
 when the actuator admits a dynamics. The open-loop ODE is exponentially stable and the open-loop infinite-dimensional system is dissipative, i.e., the energy is nonincreasing, 
but its equilibrium point is not necessarily attractive.
The feedback design is based on an extension of a finite-dimensional method, namely the \textit{forwarding method}. We propose some sufficient conditions that imply
 the well-posedness and the global asymptotic stability of the closed-loop system. 
As illustration, 
we  apply these results to a
 transport equation coupled with an ODE.  
\end{abstract}

\section{Introduction}\label{sec:intro}

This paper deals with the stabilization problem 
of systems in cascade, 
in which the first subsystem is an ordinary differential equation (ODE), and the second one is
 an infinite-dimensional system.
Another term to denote the structure under consideration is the \emph{feedforward-form}, that has been studied deeply in the context of finite-dimensional systems, 
see, e.g.,  \cite{mazenc1996adding}. 

There exist many systems in such a form.
One interesting example is the PI-controller case \cite{terrand2019adding}, where an additional finite-dimensional dynamic is added to achieve the regulation problem for an 
infinite-dimensional system.
 In our case, it is not added to obtain better performances: the finite-dimensional dynamics is imposed.
 It may represent, for instance, the dynamics of the actuator of the infinite-dimensional system 
 to be controlled. Moreover, we suppose that 
 our control is nonlinear. This situation may arise when considering saturating controllers \cite{tarbouriech2011book_saturating,marx2018stability}.
Finally, another example of application of the considered class of systems in feedforward form is the case in which an infinite-dimensional internal-model based regulator is used to solve a robust output regulation problem for a  finite-dimensional system
\cite{califano2019stability,weiss1999repetitive,paunonen2010internal,
paunonen2019stability,astolfi2021repetitive}. In this case, our proposed methodology can be employed 
to solve the stabilization problem.

Designing an explicit stabilizing feedback law in the infinite-dimensional context is not an easy task. 
To the best of our knowledge, very few is known nowadays. The backstepping method \cite{krstic_smyshlyaev_backstepping} has been proved to be applicable for many PDEs, even coupled PDE/ODE systems \cite{auriol2018delay,tang2011stabilization}. Note, however, that the method applies only on linear PDEs. It does not allow to stabilize \emph{globally} nonlinear PDEs, but only locally around the desired equilibrium \cite{cerpa_coron_backstepping,kang2018boundary}. Moreover, the backstepping method gives raise to a kernel function representing the gain of the feedback law and which solves itself a PDE that might be difficult in some cases to compute numerically. It is also worthy mentionning the construction of Lyapunov functionals with Legendre polynomials for coupled PDE/ODE systems, which leads to a hierarchy of Linear Matrix Inequalities (LMIs), and that has been applied on many systems such as the transport/ODE system \cite{safi2017tractable}, the wave/ODE system \cite{barreau2018lyapunov} or the heat/ODE system \cite{baudouin2019stability}. Note however that the Legendre method has been only applied on linear systems in 1D.

In the semigroup theory (see, e.g., \cite{tucsnak2009observation} for an introduction to the linear theory and \cite{miyadera1992nl_sg} for the nonlinear one), there exists also some design method for the stabilization of \emph{linear} systems, but it is difficult most of the time to obtain explicit gains for the feedback law. Let us mention the method provided in \cite{urquiza2005rapid}, which allows to achieve rapid stabilization of the closed-loop system and which relies on the Grammian operator. This method has been applied on the linear Korteweg-de Vries (KdV) equation \cite{cerpa2009rapid}. For coupled systems, it is also worthy mentionning some researchs related to the ouput regulation \cite{paunonen2010internal,paunonen2019stability}, where some coupled systems arise when adding the dynamics, or even \cite{feng2020actuator}, which deals with similar systems that the one we are facing with, but where the actuator dynamics is linear. 

It is important then to emphasize on the fact that, in this paper, we aim at proposing a design methods to stabilize \emph{globally} coupled systems admitting a class of \emph{input nonlinearities}. A typical example of such a nonlinearity is the saturation, but we consider a larger class of nonlinearity, namely \emph{cone-bounded nonlinearities}, as in \cite{map2017mcss}. Note that very few is known about methods for the stabilization of nonlinear PDEs, but let us mention \cite{mcpa2017siam}, which deals with the stabilization of a nonlinear KdV subject to a saturation, \cite{prieur2016wavecone,chitour2019p}, which addresses the stabilization problem of some wave equations with saturated controller; or \cite{marx2018stability,slemrod1989mcss}, which builds Lyapunov functionals for abstract systems. However, these papers do not deal with coupled systems. Let us however mention \cite{daafouz2014nonlinear}, where a coupled PDE/ODE system is considered with a saturated input.

As mentionned earlier, our method is \emph{nonlinear}. It is based on the so-called \emph{forwarding} method, which has been first introduced in the finite-dimensional context \cite{mazenc1996adding,kaliora2005stabilization,benachour2013forwarding}, and it applies on cascade systems, with particular stability properties: the first subsystem needs to be 
globally stabilizable at the origin, while the second one is dissipative. The forwarding method has been applied in  other contexts than the stabilization one, such as in the design of integral-based controllers \cite{astolfi2016integral}. To the best of our knowledge, this method has been extended in an infinite-dimensional setting only in \cite{terrand2019adding}, but, as mentioned earlier, the context is different. 

The forwarding method that we use gives raise to a gain which is a solution to an infinite-dimensional version of the well-known \emph{Sylvester equation}. We show in this paper that this gain is not difficult to obtain in some simple cases, since it can be obtained as the solution of some well-know boundary value problems. Our method needs some \emph{approximate observability} condition to apply the LaSalle's Invariance principle. This condition might seem difficult to check in general, but we also provide a sufficient condition for having such a property provided that the infinite-dimensional system under consideration is conservative (but not attractive at the origin). This condition is known in the output-regulation context as a \emph{non-resonance} condition see, e.g., \cite[Chapter 1.4]{isidori2012robust}). Let us however mention that we are not able to impose a decay rate for the trajectory of the closed-loop system, in contrast with the backstepping method. 

The paper is organized as follows. In Section \ref{sec:problem}, we introduce the system under consideration together with the functional setting which will be used all along the paper. In Section \ref{sec:forwarding}, we explain the forwarding approach and introduce the infinite-dimensional version of the Sylvester equation. In Section \ref{sec:well-posed}, we state and prove the well-posedness of the closed-loop system by using some nonlinear semigroup theory results. In Section \ref{sec:stability}, we state and prove the global asymptotic stability result. Section \ref{sec:observability} provides and proves a sufficient condition for the approximate observability in infinite for our system in cascade. Finally, Section \ref{sec:conclusion} collects some further research lines to be followed.

 \paragraph{Notation:} Set $\RR_+ = [0,\infty)$. We denote by $|\cdot|$ the Euclidean norm and by $\Vert\cdot\Vert$ the induced matrix norm. Given two Hilbert spaces $H_1$ and $H_2$, the space $\cL(H_1,H_2)$ denotes the space of functions bounded from $H_1$ to $H_2$, and $\cL(H_1) = \cL(H_1, H_1)$. Given $N\in\mathbb{N}$, for a function $w:\: (t,x)\in\RR_+\times [0,1]\mapsto w(t,x)\in\RR^N$, the notation $w_t$ (resp. $w_x$) denotes the partial derivative of $w$ with respect to the variable $t$ (resp. with respect to the variable $x$). We keep this notation
 both   for 
  the weak and the strong definition of partial derivatives.


\section{Problem Statement}\label{sec:problem}

Let $H$ be a Hilbert space equipped with a scalar product $\langle \cdot,\cdot\rangle_H$ and the corresponding norm $\Vert \cdot \Vert_H$.
In this paper we are interested in the stabilization (at the origin) problem for systems that can be described as a cascade of two systems reading as follows
\begin{equation}
\label{system}
\left\{
\begin{array}{rcl}
\frac{\d}{\d t} z &= &A z + B\sigma(u),
\\[.5em]
y&=&Cz,
\\[.5em]
\frac{\d}{\d t} w &=& S w + \Gamma y,
\\[.5em]
z(0)&=&z_0,\: w(0)=w_0,
\end{array}
\right.
\end{equation}
where $A\in \RR^{n\times n}$, $B\in \RR^{m\times n}$, $C\in \RR^{n\times p}$, $u\in \RR^m$ is the control input and  $\sigma:\RR^m\rightarrow \RR^m$ is a cone-bounded nonlinearity, that we will define in a rigorous way below. In system \eqref{system}, $S:D(S)\subseteq H\rightarrow H$ is a (possibly unbounded) operator, with $D(S)$ densely defined in $H$. Defining $H_{-1}$ as the completion of $H$ with respect to the norm $\Vert w\Vert_{-1}:=\Vert (\beta I_H -S)^{-1}w\Vert_{H}$, where $\beta$ is in the resolvent of $S$, we suppose that $\Gamma \in \cL(\RR^p,H_{-1})$ (i.e. $\Gamma$ is a bounded operator from $\RR^p$ to $H_{-1}$).

Recall that $H_{-1}$ is the dual of $D(S^*) = \{w\in H \mid \sup_{v\in D(S)\setminus\{0\}}|\langle Sv, w\rangle_H|/\|v\|_H\}$
with respect to the pivot space $H$.

In the following, we use the state space $X:=\RR^n\times H$ equipped with the  norm defined by $\Vert (z,w)\Vert_X = |z| + \Vert w\Vert_{H}$. We also consider the following set of assumptions, on which our results rely.

\begin{assumption}\label{assumption1}
The following statements hold.
\begin{enumerate}[label=(\roman*),leftmargin=*, series=hyp]
\item\label{item1} The operator $S:D(S)\subseteq H\rightarrow H$ generates a strongly continuous semigroup of contractions, that is denoted by $(\TT(t))_{t\geq 0}$. 

\item\label{item2} The matrix $A$ is Hurwitz.

\item\label{item3} The spectra of $S$ and $A$ are disjoint and nonempty.


\end{enumerate}
\end{assumption}

In system \eqref{system}, the state component $z$ lives in  $\RR^n$, while  $w$ is a (possibly infinite-dimensional) state living in the Hilbert space $H$. System \eqref{system} can be viewed as an infinite-dimensional control system in which the $z$-dynamics represent the actuator's dynamics. This class of systems may arise also in output regulation and repetitive control problems (see, e.g., \cite{paunonen2010internal,weiss1999repetitive,astolfi2021repetitive}): in these cases, the $w$-dynamics represents the  state of a dynamical feedback which is used to achieve a desired control goal for the controlled finite-dimensional $z$-subsystem.


\begin{definition}[Cone-bounded nonlinearity]
\itshape
\label{def-cb}
A continuous function $\sigma:\: \RR^m\rightarrow \RR^m$ is said to be a \textit{cone-bounded nonlinearity} if it satisfies the following properties.
\begin{itemize}
\item[1.] For all $s\in\RR^m$, $\sigma(s)=0$ if and only if $s=0$.\label{item01}
\item[2.] It is monotone\footnote{This terminology is used in functional analysis to denote monotone operators. If $m=1$, then $\sigma$ is nondecreasing.}, which means that, for every $s_1,s_2\in\RR^m$
\begin{equation}\label{sigma_monotony}
(\sigma(s_1)-\sigma(s_2))^\top (s_1-s_2)\geq 0.
\end{equation}
\item[3.] It is linearly bounded, that is: there exists a positive constant $L$ such that, for all 
$s\in \RR^m$
\begin{equation}\label{sigma_lipschitz}
|\sigma(s)|\leq L|s|.
\end{equation}
\end{itemize}
\end{definition}
The following examples are borrowed from \cite{map2017mcss}.
\begin{example}[Examples of cone-bounded nonlinearities]
\begin{itemize}
\item[1.] Any linear mapping $\sigma(s)=\mu s$, where $\mu$ is a positive constant, is a cone-bounded nonlinearity.
\item[2.] The so-called saturation function for $s\in \RR$:
\begin{equation}
\label{eq:saturation}
\sat_{\bar u}(s):=\left\{
\begin{aligned}
&-\bar u &\text{ if  }\hspace{0.2cm} &
s \leq -\bar u, \\
& s &\text{ if  }\hspace{0.2cm} &|s| \leq \bar u,\\
& \bar u  &\text{ if  }\hspace{0.2cm} &s\geq \bar u,\\
\end{aligned}
\right.
\qquad
\qquad
\forall s \; \in \RR
\end{equation}
where $\bar u$ is a positive constant, is a cone-bounded nonlinearity. By using \eqref{eq:saturation}, we can also define a saturation function when $s\in\RR^m$ as follows: 
\begin{equation}
\label{eq:saturation_m}
\sat_{\overbar U }(s):=
\begin{pmatrix}
\sat_{\bar u_1}(s_1)
\\
\vdots
\\
\sat_{\bar u_m}(s_m)
\end{pmatrix},
\qquad
\qquad
\forall \; s:= (s_1, \ldots, s_m)^\top\in \RR^m
\end{equation}
where $\overbar U := (\bar u_1, \ldots, \bar u_m)$ is a vector of positive constants $\bar u_i>0$, for all $i=1, \ldots, m$. The saturation function \eqref{eq:saturation_m} is still a cone-bounded nonlinearity.
\item[3.] The function 

\begin{equation}
\sigma(s):\: s\in\RR \mapsto \sat_{\bar{u}}(\varphi(s)), 
\end{equation}
where $\bar{u}>1$

and $\varphi$ is defined as follows 
\begin{equation}
\varphi(s):=\left\{
\begin{aligned}
& -\sqrt{|s| - 1}-1 &\text{ if } &s <-1,\\
& s &\text{ if } &|s|\leq 1, \\
&\sqrt{s-1} + 1 &\text{ if } & s > 1 
\end{aligned} 
\right.
\end{equation}
is also a cone-bounded nonlinearity, but it is not globally nor locally Lipschitz, because of the square root function in its definition. \hfill$\triangle$
\end{itemize}
\end{example}

\begin{remark}[About $A$ Hurwitz]
It is well-known that in presence of input nonlinearities, 
global stabilization of the origin of a finite-dimensional system of the form 
$$
\dfrac{\d}{\d t}z = Az + B\sigma(u)
$$ 
is not possible without any restrictive
condition on the spectrum of $A$ and the 
stabilizability property
of the pair $(A,B)$. 
For instance, if $\sigma$ is a saturation function,
a necessary condition is that 
the eigenvalues of $A$ must not have positive real parts
 \cite[Chapter 1.6.2.1, Theorem 1.2]{tarbouriech2011book_saturating}.
 In our context, due to the presence of the 
 infinite-dimensional properties of system 
 \eqref{system}, we suppose a slightly more
 stringent condition in Assumption~\ref{assumption1}-(ii), that is,   $A$ Hurwitz. An extra condition
 on $B$ will be required later on, in 
 Theorem~\ref{thm-GAS}.
Note, however, that, in absence of input nonlinearity, that is, $\sigma(u)=u$, Assumption~\ref{assumption1}-(ii) {\color{blue}can} be  replaced by supposing the pair $(A,B)$ to be stabilizable.
\end{remark}

We provide now an example of a transport equation coupled with an ODE which corresponds to our functional setting. Such an example will be also used in the following sections as an illustration of the proposed design.

\begin{example}[Transport equations coupled with an ODE]
\label{ExPart1}
Consider
\begin{equation}
\label{hyperbolic}
\left\{
\begin{array}{ll}
\dot{z}(t) = Az(t) + B\sigma(u(t)), & t\in \RR_+\\
w_t(t,x) + \Lambda w_x(t,x) = 0,  & (t,x)\in\RR_+\times [0,1]\\
w^+(t,0) = D_0 w^-(t,0), & t\in \RR_+,\\
w^-(t,1) = D_1 w^+(t,1) + Cz(t),\quad & t\in \RR_+\\
z(0)= z_0,\: w(0,x) = w_0(x), & x\in [0,1],
\end{array}
\right.
\end{equation}
where $\Lambda=\diag(\lambda_1,\ldots,\lambda_N)$ is such that $\lambda_i>0$ (resp. $\lambda_i<0$), for every $i\in\lbrace 1,\ldots,k\rbrace$ (resp. $\lambda_i>0$ for every $i\in\left\lbrace k+1,\ldots N\right\rbrace$), $w^+(t,x)$ (resp. $w^-(t,x)$) corresponds to the $k$ first components of $w$ (resp. to the $N-k$ last components of $w$), $C$, $D_0$ and $D_1$ are matrices of appropriate dimension. 
The matrix $A$ is Hurwitz and the function $\sigma$ is supposed to be a cone-bounded nonlinearity.
The boundary control problem \eqref{hyperbolic} can be rephrased as \eqref{system} according to \cite[Chapter 10]{tucsnak2009observation}.
Following the steps of \cite[Example 10.1.9]{tucsnak2009observation}, operators $S$ and $\Gamma$ can be determined as follows.
The state space is chosen to be $X:=\RR^n\times H$, $H=L^2(0,1;\RR^N)$. The operator $S$ is given by $Sw:=-\Lambda w^\prime$
with domain $D(S):=\lbrace w\in H^1(0,1;\RR^N)\mid w^+(0) = D_0w^{-}(0),\: w^{-}(1)=D_1 w^{+}(1)\rbrace$,
and the operator $\Gamma$ is the delta function at $x=1$ in $\cL(\RR^p,H_{-1})$, \emph{i.e.}, $\langle v, \Gamma y \rangle_{D(S^*), H_{-1}} = v(1) y$ for all $y\in\RR^p$ and $v\in D(S^*)$ where $\langle \cdot, \cdot \rangle_{D(S^*), H_{-1}}$ is the dual product.

Now, under some assumptions on the matrices $D_0$ and $D_1$ to be given later on, let us prove that Asssumption~\ref{assumption1}~\ref{item1} is satisfied. First, by \cite[Theorem 3.1.]{russell1978controllability}, the operator $S$ generates a strongly continuous semigroup (for all $D_0$ and $D_1$). Then, $S$ is dissipative if and only if $S$ generates a strongly continuous semigroup of contractions.
To do so, let us introduce the following real inner product on $H$: for all $w_1, w_2\in H$,
$$
\langle w_1,w_2\rangle_{H}:=\int_0^1 w_1(x)^\top \tilde{\Lambda} w_2(x) dx,\quad \tilde{\Lambda}:= \diag\left(\frac{1}{|\lambda_1|},\dots,\frac{1}{|\lambda_N|}\right).
$$
Showing that $S$ is dissipative reduces to prove that $\langle S w,w\rangle_H \leq 0$. We have
\begin{equation*}
\langle Sw,w\rangle_H = - \int_0^1 \big((w^{+})^\prime (x)\big)^\top w^{+}(x) dx + \int_0^1 \big((w^{-})^\prime (x)\big)^\top w^{-}(x) dx,
\end{equation*}
One deduces by integration by parts that
\begin{align*}
\langle Sw,w\rangle_H = &~ \frac{1}{2}\left (w^{-}(0)^\top D_0^\top D_0 w^{-}(0) - w^-(0)^\top w^-(0)\right) \\
&+ \frac{1}{2}\left( w^{+}(1)^\top D_1^\top D_1 w^{+}(1) - w^+(1)^\top w^+(1)\right).
\end{align*}
Hence, if matrices $D_0$ and $D_1$ satisfy $D_0^\top D_0 - I_{M} \preceq 0$ and $D_1^\top D_1 - I_{N-M} \preceq 0$, one can deduce that $S$ is dissipative.
Conversely, if these matrix inequalities are not satisfied, one can pick $w^-(0)$ and $w^+(1)$ such that $\langle Sw,w\rangle_H\geq0$, hence $S$ does not generate a strongly continuous semigroup of contractions.


\hfill$\triangle$
\end{example}

Our objective in this paper is to propose a design procedure to stabilize \eqref{system}. We follow a \emph{forwarding} strategy, that has been applied in the finite-dimensional context in many situations \cite{mazenc1996adding}, and that has been extended recently to some infinite-dimensional systems in a different context \cite{terrand2019adding}.

\section{The Forwarding Approach}\label{sec:forwarding}

Inspired by the forwarding strategy and in particular by the change of coordinate approach proposed
in \cite[Section IV]{mazenc1996adding}, we look for an operator $M:\RR^n\rightarrow H_{-1}$ satisfying
\begin{equation}
\label{eq-Sylvester}
SM - MA = -\Gamma C,
\end{equation}
where $S$ is understood as the extension of the operator $S$ in $H_{-1}$.
Since the spectra of $A$ and $S$ are disjoint 
(see Assumption~\ref{assumption1}~\ref{item3}),
we conclude 
that there exists  a unique $M:\RR^n\rightarrow H_{-1}$ solution of \eqref{eq-Sylvester}, see 
 \cite[Lemma 22]{phong1991operator}.
Also, $M$ takes values in $H$
as proved
in \cite[Lemma 4.2.]{feng2020actuator},
and, in particular, 
we have $SMz+\Gamma Cz\in H$
for all $z\in\RR^n$.

Let us recall the proof briefly. Take 
$\mu$ in the resolvent of $S$ and consider the Sylvester equation for all $z\in\RR^n$
$$
SMz - MAz -\mu Mz + Mz\mu = -\Gamma Cz.
$$
By a simple computation, one can prove that, for all
 $z\in\RR^n$,
\begin{equation}\label{sol:sylv}
    Mz = (\mu I_H - S)^{-1} M (\mu I_{n} - A)z + (\mu I_H-S)^{-1} \Gamma Cz.
\end{equation}
Since $(\mu I_H - S)^{-1}\in\mathcal{L}(H_{-1},H)$, one can deduce from the above expression that $Mz\in H$ for all $z\in\RR^n$.


\begin{example}[Sylvester equation for \eqref{hyperbolic}]
In Example~\ref{ExPart1}, the Sylvester equation~\eqref{eq-Sylvester} can be rewritten as follows:
\begin{equation*}\label{eq:SylvesterTwoBoundary}
\left\{
\begin{array}{ll}
-\Lambda M^\prime(x) - M(x)A = 0, &\quad  x\in [0,1]\\
M^+(0) =D_0 M^-(0)\\
M^{-}(1) =  D_1 M^+(1) + C
\end{array}
\right.
\end{equation*}
The solution $M: [0,1]\to \RR^{N\times m}$ to this two-boundary value problem can be computed explicitly using some vectorization operators.
Consider the case of the one-dimensional transport equation driven by a scalar ODE:
\begin{equation}
\label{scalar-hyperbolic}
\left\{
\begin{array}{ll}
\dot{z} = -a z + \sigma(u),&\: t\in\mathbb{R}_+\\
w_t(t,x) + \lambda w_x(t,x) = 0,&\: (t,x)\in\mathbb{R}_+\times [0,1]\\
w(t,0) = w(t,1) + c z,&\: t\in\mathbb{R}_+\\
z(0) = z_0,\: w(0,x) = w_0(x),&\: x\in [0,1],
\end{array}
\right.
\end{equation}
with $a$, $c$ and $\lambda$ three positive constants.

System~\eqref{scalar-hyperbolic} can be seen as a subcase of \eqref{hyperbolic} with $N=2$ and $n=m=p=1$. Indeed, taking $\lambda = 2\lambda_1=-2\lambda_2$, $D_0=D_1=1$, $A=-a$, $B=1$, $C=c$ and $$
w_0(x) = \left\{\begin{array}{ll}w^+_0(2x)  \quad & \text{ if } x<\frac{1}{2} \ ,
\\w^-_0(2(1-x)) \quad & \text{ if } x>\frac{1}{2} \ ,
\end{array}
\right.
$$
any $(z, (w^+, w^-))$ satisfies \eqref{hyperbolic} if and only if $(z, w)$ is a solution of \eqref{scalar-hyperbolic}, where 
$$
w(t,x) = \left\{ \begin{array}{ll}w^+(t,2x) \quad &\text{ if } x<\frac{1}{2} \ ,\\
w^-(t,2(1-x)) \quad &\text{ if } x>\frac{1}{2} \ .
\end{array}
\right.
$$

Then $M$ is the unique solution of
\begin{equation}
\label{sylvester-scalar}
\left\{
\begin{array}{ll}
\lambda M^\prime(x) = a M(x),&\quad x\in [0,1]\\
M(0) = M(1) + c,
\end{array}
\right. 
\end{equation}
i.e., for all $x\in[0, 1]$,
$$
M(x) =\frac{c}{1-\exp\left(\frac{a}{\lambda}\right)} \exp\left(\frac{a}{\lambda} x\right),
$$
\hfill$\triangle$
\end{example}

Now go back to the general case \eqref{system}.The fact that $A$ is Hurwitz implies that there exists a symmetric positive definite matrix $P$ satisfying
\begin{equation}
\label{lyapunov-inequality}
PA + A^\top P = -I_{\RR^n}.
\end{equation}
We denote by $p_{\min}>0$ (resp. $p_{max}>0$) its smallest (resp. largest) eigenvalue.
Following the forwarding approach, let us first 
introduce the following candidate Lyapunov functional
\begin{equation}
\label{lyapunov-function}
V(z,w):= z^\top P z + \Vert w-Mz\Vert_H^2,
\end{equation}
where $M$ is the operator previously defined in \eqref{eq-Sylvester}. Then $V$ induces a scalar product on $X$ defined by

\begin{equation}
\label{lyapunov-scalar-product}
\left\langle \begin{bmatrix}
z_1 & w_1
\end{bmatrix}^\top,\begin{bmatrix}
z_2 & w_2
\end{bmatrix}^\top \right\rangle_V := z_1^\top P z_2 + \langle w_1-Mz_1,w_2-Mz_2\rangle_H.
\end{equation}

The norm induced by $V$, {\color{blue} denoted in the following by 
$\|\cdot\|_V$,}
is equivalent to $\|\cdot\|_X$. Indeed, for all $(z, w)\in X$, we have:
\begin{align*}
    V(z, w)
    &=
    z^\top P z + \Vert w\Vert_H^2 + \Vert Mz\Vert_H^2 - 2 \langle w, Mz\rangle\\
    &\leq p_{\max} |z|^2
    + 2(\|w\|^2_H + \|Mz\|^2_H)\\
    &\leq \max\left(2, p_{\max} +2\|M\|^2_{\cL(\RR^n, H)}\right)\|(z, w)\|^2_X
\end{align*}
and

\begin{align*}
    V(z, w)
    &=
    z^\top P z + \Vert w\Vert_H^2 + \Vert Mz\Vert_H^2 - 2 \langle w, Mz\rangle\\
    &\geq p_{\min} |z|^2
    + (1 - \varepsilon)\Vert w\Vert_H^2
    + \left(1 - \frac{1}{\varepsilon}\right)\Vert Mz\Vert_H^2\\
    &\geq \left(p_{\min}+\left(1-\frac{1}{\varepsilon}\right)\|M\|^2_{\cL(\RR^n, H)}\right)
    |z|^2 + (1-\varepsilon)\Vert w\Vert_H^2\\
    &\geq \min\left(\frac{p_{\min}}{2}, \frac{p_{\min}}{p_{\min} + 2\|M\|^2_{\cL(\RR^n, H)}}\right) \|(z, w)\|^2_X
\end{align*}
by choosing 
$$
\varepsilon = \frac{2\|M\|^2_{\cL(\RR^n, H)}}{p_{\min}+2\|M\|^2_{\cL(\RR^n, H)}} \ .
$$


Let $(z, w)$ be a sufficiently regular solution to \eqref{system}, i.e., $(z,w)\in D(\cA)$, where
\begin{equation}\label{domA}
D(\cA): = \lbrace (z,w)\in X\mid Sw + \Gamma C z \in H\rbrace.
\end{equation}
Since $(\TT(t))_{t\geq 0}$ is a strongly continuous semigroup of contractions, $S$ is a dissipative operator.
Then,
\begin{align}
\frac{\d}{\d t} V(z,w) = &  ~2z^\top P(A z +
  B \sigma(u)) 
  \nonumber\\
&+ \langle Sw-MAz + \Gamma Cz,w-Mz\rangle_H \nonumber
\\
& + \langle w-Mz,Sw-MAz + \Gamma Cz\rangle_H \nonumber
\\
& - \langle M B\sigma(u),w-Mz\rangle_H
- \langle w-Mz,M B\sigma(u)\rangle_H,
\nonumber\\
= &- z^\top z + 2z^\top P B \sigma(u)
\tag{by \eqref{eq-Sylvester} and \eqref{lyapunov-inequality}}\\
& + \langle Sw-SMz,w-Mz\rangle_H + 
\langle w-Mz, Sw-SMz\rangle_H  \nonumber
\\
& - \langle M B\sigma(u),w-Mz\rangle_H
- \langle w-Mz,M B\sigma(u)\rangle_H,
\nonumber\\
\leq & -z^\top z +
 2z^\top P B \sigma(u)
 \tag{since $S$ is dissipative}\\
 &- \langle M B\sigma(u),w-Mz\rangle_H - \langle w-Mz,M B\sigma(u)\rangle_H,
\nonumber\\
\leq &-z^\top z + 2 \big[z^\top P - (M^*(w-Mz))^\top
\big] B\sigma(u).
\label{derivative-lyap}
\end{align}
Hence, a natural candidate control $u$ is:
\begin{equation}
\label{feedback}
u=-B^\top[Pz - M^*(w-Mz)]
\end{equation}
yielding
$$
\frac{\d}{\d t} V(z,w) \leq - z^\top z- 2 u^\top \sigma(u) \leq 0.
$$ 
The previous inequality being not strict, i.e., we don't have any negative term in $\|(z,w)\|_X$, we cannot conclude any attractivity result without a  more precise analysis. This will be done in Section~\ref{sec:stability}.
The closed-loop system \eqref{system}-\eqref{feedback} given by
\begin{equation}
\label{cl-system}
\left\{
\begin{aligned}
&\frac{\d}{\d t} z = A z + B\sigma\left(-B^\top \left(Pz - M^*(w-Mz)\right)\right),\\
&\frac{\d}{\d t} w = S w + \Gamma C z,\\
&z(0)=z_0,\: w(0)=w_0,
\end{aligned}
\right.
\end{equation}
can be rewritten as
\begin{equation}
\label{operator-A}
\left\{
\begin{aligned}
&\frac{\d}{\d t} \zeta = \cA(\zeta),\\
&\zeta(0) = \zeta_0,
\end{aligned}
\right.
\end{equation}
where $\zeta=\begin{bmatrix}
z & w
\end{bmatrix}^\top$ and

$$
\cA(\zeta) = \begin{bmatrix}
Az + B\sigma\left(-B^\top\left( Pz - M^*(w-Mz)\right)\right)  \\
Sw + \Gamma Cz
\end{bmatrix}.
$$

The domain of $\cA$ is $D(\cA)$ (given in \eqref{domA})
equipped with the graph norm. In the next section, we first study the 
existence and uniqueness of solutions of \eqref{cl-system}.

\begin{example}[Control design for \eqref{scalar-hyperbolic}]
For system \eqref{scalar-hyperbolic}, the control law \eqref{feedback} can be explicitly design by using the function $M$ defined in \eqref{sylvester-scalar}. In particular, in this case, $M^*$ is given by 
$M^*:\: L^2(0,1)\ni w\mapsto \int_0^1 M(x) w(x) \d x\in\RR$. Hence, the feedback \eqref{feedback} reads, for all $t\geq 0$
$$
u(t) = - \dfrac{1}{2a}z+\int_0^1 M(x) [w(t,x)- M(x)z]\d x.
$$
See also \cite{marx2020forwarding}.
\hfill$\triangle$
\end{example}

\begin{remark}[Nonlinear design and small control]
Note that 
the feedback law \eqref{feedback}
can be also modified in 
$$
u = \psi\Big(- B^\top Pz - M^*(w-Mz)\Big),
$$
with $\psi$ being any desired function so that 
 the  composition $\sigma\circ \psi$ is still a 
cone-bounded
function satisfying Definition~\ref{def-cb}.
Indeed, the resulting closed-loop system 
would still read as \eqref{cl-system} and all the 
forthcoming results  apply.
For instance, if $u$ takes value in $\RR$, 
it suffices to select $\psi$ as any 
cone-bounded function
satisfying Definition~\ref{def-cb}.
 Recall indeed that  the composition of two
monotonic 
 functions $f,g:\RR^n\to \RR$, is still a 
 monotonic function. The other two properties
 in Definition~\ref{def-cb}
 are trivially verified.
Such a choice 
can be of   particular 
interest for the 
design of  saturated feedback laws 
with small magnitude (that is, small energy), 
see, e.g., 
\cite{tarbouriech2011book_saturating}.
Indeed, by  selecting $\psi$ as a saturation function 
of the form
\eqref{eq:saturation_m}, 
it turns out that 
the saturation level $\bar U$ can be chosen 
arbitrarily small to reduce the efforts on the
actuators (although, as a drawback, 
this may also sensitively reduce
the converge rate).
Such a feature,  highly desirable from a practical point of view, is a typical feature of forwarding approach for finite-dimensional systems (see, e.g., \cite{kaliora2005stabilization}) which is still preserved in our context, and that is 
very hard to retain when considering pole-placement techniques (e.g., \cite{feng2020actuator,paunonen2010internal})
or backstepping design
(see, e.g., \cite{krstic_smyshlyaev_backstepping,auriol2018delay,
tang2011stabilization}) in the infinite-dimensional framework.
\end{remark}

\section{Well-posedness Result}\label{sec:well-posed}

We are now in position to state our first result, which deals with the well-posedness of system \eqref{cl-system}. 

\begin{theorem}[Well-posedness and Lyapunov stability of (\ref{cl-system})]
\label{thm-wp}
Suppose Assumption~\ref{assumption1} is satisfied. Then, the following statements hold.
\begin{itemize}
\item[1.] For every initial conditions $(z_0,w_0)\in X$, there exists a unique \emph{weak solution} $(z,w)\in C^0(\RR_+;X)$ to \eqref{cl-system}. Moreover, for all $t\geq 0$,
\begin{equation}
\label{lyapunov-X}
\Vert (z(t),w(t))\Vert_V \leq \Vert (z_0,w_0)\Vert_V. 
\end{equation}
\item[2.] For every initial conditions $(z_0,w_0)\in D(\cA)$, there exists a unique \emph{strong solution} $(z,w)\in C^1(\RR_+;X)\cap C^0(\RR_+;D(\cA))$ to \eqref{cl-system}. Moreover, for all $t\geq 0$,
\begin{equation}
\label{lyapunov-D(A)}
\Vert (z(t),w(t))\Vert_{V}^2 + \Vert \cA (z(t),w(t))\Vert_{V}^2 \leq  \Vert (z_0,w_0)\Vert_{V}^2 + \Vert \cA (z_0,w_0)\Vert_{V}^2 \ .
\end{equation}
\end{itemize} 
\end{theorem}

If one proves that the operator $\cA$ defined in \eqref{operator-A} is a $m$-dissipative operator on $(X, \|\cdot\|_X)$,  according to the definition below, one can apply the result provided by \cite[Corollary 3.7, Theorem 4.20]{miyadera1992nl_sg}, and conclude that the statements of Theorem \ref{thm-wp} hold.
\begin{definition}[$m$-dissipative operators]
\label{def_dissipative}
\itshape
An operator $\cA:\: D(\cA)\subset X\rightarrow X$ is said to be  $m$-dissipative if and only if
\begin{itemize}
\item[$\bullet$] The operator $\cA$ is \emph{dissipative}, i.e.
\begin{equation}
\langle \cA(\zeta_1)-\cA(\zeta_2),\zeta_1-\zeta_2\rangle_V\leq 0,\qquad \forall \zeta_1,\zeta_2\in D(\cA).
\end{equation}
\item[$\bullet$] The operator $\cA$ is maximal, i.e. there exists $\lambda_0>0$ such that (equivalently, for all $\lambda_0>0$)
\begin{equation}
 \Ran(\lambda_0 I_X - \cA)=X,
\end{equation}
where $I_X$ is the identity operator over the Hilbert space $X$, and $\Ran$ is the range operator.
\end{itemize}
\end{definition} 

\begin{proof}[Theorem \ref{thm-wp}]
The proof of Theorem \ref{thm-wp} relies on the so-called semigroup theory (see, e.g., \cite{pazy1983semigroups} or \cite{tucsnak2009observation} for an introduction to this theory in the linear case, and \cite{miyadera1992nl_sg} in the nonlinear case). The proof is divided in two steps: we first prove that $\cA$ is dissipative and second that $\cA$ is maximal. 

\textbf{First step: $\cA$ is dissipative.}
To prove this result, we use the scalar product introduced in \eqref{lyapunov-scalar-product}, that is equivalent to the standard scalar product used for $X$.
Let $\zeta_1 = \begin{bmatrix}
z_1 & w_1
\end{bmatrix}^\top$ and $\zeta_2:=\begin{bmatrix}
z_2 & w_2
\end{bmatrix}^\top$ be in $ D(\cA)$.
We introduce $\tilde{z}=z_1-z_2$ and $\tilde{w}=w_1-w_2$,  $u_i= - B^\top(Pz_i - M^*(w_i-Mz_i))$.

Since $(\tilde{z},\tilde{w})\in D(\cA)$ and because $S\tilde{w} + \Gamma C \tilde{z}\in H$, one can compute $\langle \cA(\zeta_1) - \cA(\zeta_2),\zeta_1-\zeta_2\rangle_V$. We obtain:
\begin{align}
&\langle \cA(\zeta_1)-\cA(\zeta_2),\zeta_1-\zeta_2\rangle_V
\\ 
&\qquad\qquad  = \tilde{z}^\top A^\top P \tilde{z}
+ (\sigma(u_1) - \sigma(u_2))^\top B^\top P \tilde{z} \nonumber
\\
&\qquad\qquad \phantom{=\;} +  \langle S\tilde{w}-(MA-\Gamma C)\tilde{z} - MB(\sigma(u_1)-\sigma(u_2)),\tilde{w}-M\tilde{z}\rangle_H\nonumber
\\
&\qquad\qquad   = -\frac{1}{2}|\tilde{z}|^2 + (\sigma(u_1) - \sigma(u_2))^\top B^\top P \tilde{z}\tag{by \eqref{eq-Sylvester}}
\\
&\qquad\qquad \phantom{=\;}  +\langle S(\tilde{w} - M\tilde{z}) -MB(\sigma(u_1)-\sigma(u_2)),\tilde{w}-M\tilde{z}\rangle_H\nonumber
\\
&\qquad\qquad \leq -\frac{1}{2}|\tilde{z}|^2 + (\sigma(u_1) - \sigma(u_2))^\top B^\top P \tilde{z}\tag{$S$ is dissipative}
\\
&\qquad\qquad \phantom{=\;} -(\sigma(u_1)-\sigma(u_2))^\top B^\top M^*(\tilde{w}-M\tilde{z})\nonumber
\\
&\qquad\qquad \phantom{=\;}   -\frac{1}{2}|\tilde{z}|^2 -(\sigma(u_1)-\sigma(u_2))^\top (u_1-u_2)\nonumber
\\
&\qquad\qquad \leq 0.
\end{align}
where in the last step we used the monotony of $\sigma$ given in \eqref{sigma_monotony}.
Thus $\cA$ is a dissipative operator, concluding the first step of the proof.

\textbf{Second step: $\cA$ is a maximal operator.} Proving that $\cA$ is maximal reduces to show that, given a positive constant $\lambda$, for all $\zeta\in X$, there exists $\tilde{\zeta}\in D(\cA)$ such that
\begin{equation}
(\lambda I_X-\cA) \tilde{\zeta} = \zeta.
\end{equation}
The constant $\lambda$ will be selected later on. Let $(z,w)\in X$. We seek $(\tilde{z},\tilde{w})\in D(\cA)$ such that
\begin{equation}
\label{eq-maximal}
\left\{
\begin{aligned}
&\lambda \tilde{z} -A\tilde{z} - B\sigma\left(B^\top (-P\tilde{z} + M^*(\tilde{w}-M\tilde{z})\right) = z\\
&\lambda \tilde{w} - S\tilde{w} - \Gamma C\tilde{z} = w.
\end{aligned}
\right.
\end{equation}
Since $S$ is a dissipative operator, the real parts of its eigenvalues are non-positive, and thus any positive constant $\lambda$ is in the resolvent of $S$. Then, the operator $(\lambda I_H-S)$ is invertible. We have
\begin{equation}
\tilde{w} = (\lambda I_H-S)^{-1} [w + \Gamma C\tilde{z}]. 
\end{equation}
We can therefore rewrite the first line of \eqref{eq-maximal} as an equation depending only on $\tilde{z}$ and the given data $w$ and $z$.
\begin{equation}
\lambda \tilde{z} - A\tilde{z} - B\sigma\left(B^\top(-P\tilde{z} + M^*((\lambda I_H-S)^{-1} [w + \Gamma C\tilde{z}] - M\tilde{z}))\right) = z. 
\end{equation}
Since $A$ is Hurwitz, one may rewrite the latter equation as follows

\begin{equation}
\label{eq-max-fp}
\tilde{z} = (\lambda I_n-A)^{-1} \left[B\sigma\left(B^\top(-P\tilde{z} + M^*((I_H-S)^{-1} [w + \Gamma C\tilde{z}] - M\tilde{z}))\right) + z\right].
\end{equation}

We prove that there exists a solution to this equation by using a fixed-point strategy. Let us define the operator
\begin{equation}
\begin{aligned}
F : \RR^n &\rightarrow \RR^n\\
\bar{z} &\mapsto F(\bar{z}),
\end{aligned}
\end{equation}
where $F$ is defined by

$$
\begin{array}{rcl}
F(\bar{z})&:=& (\lambda I_n-A)^{-1} \Big[z + 
\\
&& \qquad \qquad
B\sigma\left(B^\top(-P\bar{z} + M^*((\lambda I_H-S)^{-1} [w + \Gamma C \bar{z}] - M\bar z))\right) 
\Big]
 \\
&=&(\lambda I_n-A)^{-1} \Big[ 
z + 
\\
&& \qquad \qquad 
B\sigma\big(B^\top \big(-P\bar{z} + M^*[(\lambda I_H-S)^{-1}(w - M(\lambda I_n-A)\bar z)]\big)\big) 
\Big],
\end{array}
$$
by \eqref{sol:sylv}.

Our aim is to prove that this operator admits at least one fixed-point, which will prove that there exists a solution to \eqref{eq-max-fp}. To this end, we apply the Schauder fixed-point theorem \cite[Theorem B.19]{bible_coron}. 

Consider the ball $\cB_R\subset X$ centered at $0$ of radius $R>0$. If for every $\bar{z}\in\cB_R$, we succeed to prove that $F(\bar{z})\in \cB_R$, then the Schauder fixed-point theorem applies (since $\cB_R$ is convex and compact), and the proof is achieved.
We have

\begin{equation}
\begin{aligned}
|F(\bar{z})|\leq & \Vert (\lambda I_{n}-A)^{-1}\Vert\Big[L\Vert B\Vert^2 \left(\Vert P\Vert + \Vert M^*(\lambda I_H-S)^{-1}M(\lambda I_n-A)\Vert\right) \Big]|\bar{z}| \\
&+ \Vert (\lambda I_{n}-A)^{-1}\Vert L\Vert B\Vert^2 |M^*(\lambda I_H-S)^{-1}w| +  \Vert (\lambda I_{n}-A)^{-1}\Vert |z|.
\end{aligned}
\end{equation}
Now choose $\lambda$ sufficiently large to apply the Schauder fixed-point theorem.
If $\lambda>\|A\|$, then
$\|(\lambda I_{n}-A)^{-1}\| \leq \frac{1}{\lambda-\Vert A\Vert}$ (see, \emph{e.g.}, \cite[Lemma 2.2.6]{tucsnak2009observation}).
Hence
$\|(\lambda I_{n}-A)^{-1}\|\to0$ as $\lambda\to+\infty$.
Moreover,
$\Vert (\lambda I_H-S)^{-1}M(\lambda I_n-A)\Vert$ remains bounded.
Indeed, $\Vert(\lambda I_n-A)\Vert \leq \lambda + \|A\|$
and $\Vert (\lambda I_H-S)^{-1}\Vert \leq \frac{m}{\lambda-\omega}$
by \cite[Corollary 2.3.3]{tucsnak2009observation} for some positive constants $m, \omega$ and all $\lambda>\omega$.
We can therefore select $\lambda$ sufficiently large such that
$$
\Vert (\lambda I_{n}-A)^{-1}\Vert\Big[L\Vert B\Vert^2 \left(\Vert P\Vert + \Vert M^*(\lambda I_H-S)^{-1}M(\lambda I_n-A)\Vert\right) \Big] \leq \frac{1}{2}
$$
Then, if one selects $R$ satisfying
\begin{equation*}
\Vert (\lambda I_{n}-A)^{-1}\Vert L\Vert B\Vert^2 |M^*(\lambda I_H-S)^{-1}w| +  \Vert (\lambda I_{n}-A)^{-1}\Vert |z|
\leq\frac{R}{2}
\end{equation*}
we obtain  $|F(\bar{z})|\leq R$.

Applying Schauder fixed-point theorem \cite[Theorem B.19.]{bible_coron}, the operator $F$ admits a fixed-point, which implies that the operator $\cA$ is maximal.
Since  $\cA$ is dissipative and maximal, we obtain that $\cA$ is $m$-dissipative, which concludes the proof of Theorem \ref{thm-wp}.\qed
\end{proof}

\section{Global Asymptotic Stability Result}\label{sec:stability}


To prove the global asymptotic stability of the origin of \eqref{cl-system}, we will invoke an infinite-dimensional version of the LaSalle's Invariance Principle \cite[Theorem 3.1.]{slemrod1989mcss}. But such a result is not sufficient to deduce the global asymptotic stability of the closed-loop system. As in the finite-dimensional context, we need an observability condition for \eqref{cl-system}. To this end, we recall below the definition of infinite-time approximate observability given in \cite[Definition 6.5.1.]{tucsnak2009observation}.

\begin{definition}[Infinite-time approximate observability]
\label{def-approximate}
\itshape
Let $H$ and $Y$ be two real Hilbert spaces,
$S:\: D(S)\subset H\rightarrow H$ be the generator of a strongly continuous semigroup $(\TT(t))_{t\geq 0}$
and $\cC\in\cL(X,Y)$.
Let $\mathbf{\psi}\in \cL(D(S),L^2([0,\infty);Y))$ be defined by
\begin{equation}\label{def-psi}
    (\mathbf{\psi} \zeta_0)(t) = \cC\TT(t) \zeta_0,\quad \forall\zeta_0\in D(S)
\end{equation}
for all $t\geq0$.
Then, the pair $(S,\cC)$ is approximately observable in infinite-time if and only if $\Ker \mathbf{\psi} = \lbrace 0\rbrace$.
\end{definition} 

\begin{theorem}[Global asymptotic stability]
\label{thm-GAS}
Suppose that Assumption~\ref{assumption1} holds.
Moreover, assume the following:\\[-1.5em]
\begin{enumerate}[label=(\roman*),leftmargin=*, resume=hyp]
      \item\label{item-resolvant} $S$ has compact resolvent;
        \item\label{pair-observable}
    $(S,B^\top M^*)$ is approximately observable in infinite-time.
\end{enumerate} Then, the origin of \eqref{cl-system} is globally asymptotically stable in $X$, that is, for every initial condition $(z_0,w_0)\in X$, the origin is Lyapunov stable and $$\lim_{t\rightarrow +\infty}\Vert (z(t),w(t))\Vert_X=0.$$
\end{theorem}


\begin{remark}[Compact resolvent]
\label{rem-compact}
According to \cite[Proposition 4.24]{cheverry2019handbook}, it is worthy mentionning that there is an equivalence between the compactness of the resolvent and the precompactness of the positive orbit, a property required to apply LaSalle's Invariance Principle \cite[Theorem 3]{slemrod1989mcss}.
\end{remark}

\begin{remark}
According to \cite[Theorem 4.1.5]{curtain1995introduction}, since $S$ is unbounded and $B^\top M^*$ takes values in a finite-dimensional space, the pair $(S,B^\top M^*)$ cannot be exactly observable.
For this reason, we rely on an approximate observability hypothesis.
\end{remark}
\begin{remark}
\label{rem-pair}
Checking whether \ref{pair-observable} holds might be difficult in practice, except in some simple cases, as the one provided in Example \ref{example-observability-scalar}. Unfortunately, this condition introduces the operator $M$, defined in \eqref{eq-Sylvester}, and which does not depend directly on the data of the initial problem given in \eqref{system}. Nevertheless, in the case where $S$ is skew-adjoint with compact resolvent, we provide a necessary condition which does not depend on $M$, see below in Proposition~\ref{thm:observability}.
\end{remark}

\begin{example}[Checking \ref{pair-observable} for \eqref{scalar-hyperbolic}]
\label{example-observability-scalar}
We aim at proving that property \ref{pair-observable} holds for the system \eqref{scalar-hyperbolic} introduced in Example~3.
Let $\zeta_0 = (z_0, w_0)\in \Ker \psi$ where $\psi$ is given by \eqref{def-psi} and $\cC = B^\top M^*$.
In this case, $M^*$ is given by $M^*:\: L^2(0,1)\ni w\mapsto \int_0^1 M(x) w(x) \d x\in\mathbb{R}$. Recalling that $B^\top=1$, the equation $B^\top M^* \TT(t) w_0(x)=0$ for all  $t\geq 0$ reduces to
\begin{equation}
\label{scalar-observability}
\int_0^1 M(x) w(t,x) \d x=0,
\end{equation}
where $w$ is the solution to
\begin{equation}
\label{scalar-hyperbolic-omega}
\left\{
\begin{array}{ll}
w_t + \lambda w_x = 0,\quad&(t,x)\in\mathbb{R}_+\times [0,1],\\
w(t,0) = w(t,1),\quad &t\in\mathbb{R}_+\\
w(0,x) = w_0(x),\quad & x\in [0,1].
\end{array}
\right.
\end{equation}
Using \eqref{sylvester-scalar}, one obtains therefore

\begin{equation}
\int_0^1 M(x) w(t,x) \dx = \frac{\lambda}{a} \int_0^1 M^\prime(x) w(t,x) \dx.
\end{equation}

Performing an integration by parts leads to

\begin{multline*}
\frac{\lambda}{a} \int_0^1 M^\prime(x) w(t,x) \dx = \frac{\lambda}{a} \Big[[M(1)w(t,1) - M(0)w(t,0)] \\ -\int_0^1 M(x) w_x(t,x) \dx\Big]\ .
\end{multline*}
By noting that \eqref{scalar-observability} implies  $\int_0^1 M(x) w_t(t,x) \dx = 0$, we obtain $$\int_0^1 M(x) w_x(t,x) \dx = 0.$$ 
Using the boundary conditions of \eqref{scalar-hyperbolic-omega} and \eqref{sylvester-scalar}, one can deduce
that $w(t,1)=0$. Then, following the proof of \cite[Theorem 2]{marx2020forwarding}, one can prove that $w(t,x) = 0$, which shows that \ref{pair-observable} holds for \eqref{scalar-hyperbolic}. Moreover, from standard Sobolev injections results, the canonical embedding from $D(S)$ into $H$ is compact, i.e., $S$ has compact resolvent. Then, Theorem \ref{thm-GAS} applies, and $(0,0)$ in $\mathbb{R}^n\times H$ is globally asymptotically stable in the $\mathbb{R}^n\times H$-topology. \hfill$\triangle$
\end{example}


\begin{proof}[Theorem \ref{thm-GAS}] We prove the statement of the theorem for initial conditions $(z_0, w_0)$ in $D(\cA)$. The result follows for all initial conditions in $X$ by a standard density argument (see e.g. \cite[Lemma 1]{map2017mcss})).
Consider the Lyapunov function defined in \eqref{lyapunov-function}. Denoting  
$u=-B^\top (Pz - M^*(w-Mz))$,
 and performing the same computations than in \eqref{derivative-lyap}, with the control defined in \eqref{feedback}, one therefore obtains 
\begin{equation}
\frac{\d V}{\d t}(z,w) \leq -z^\top z - 2u^\top \sigma(u)
\end{equation}
for every $(z_0,w_0)\in D(\cA)$.
Then, for all $t\geq 0$,
\begin{equation}
V(z(t),w(t)) - V(z_0,w_0) \leq -\int_0^t z(s)^\top z(s) \d s - 2\int_0^t u(s)^\top \sigma(u(s)) \d s.
\end{equation}
This implies that
\begin{equation}
\int_0^\infty z(s)^\top z(s) \d s <+\infty,\quad \int_0^\infty u(s)^\top \sigma(u(s)) \d s<+\infty.
\end{equation}
Since the trajectories of \eqref{cl-system} are bounded in $D(\cA)$ \eqref{operator-A}, and $\sigma$ is linearly bounded, then $z$ and 
$$
s\mapsto \left(z(s)^\top P-M^* (w(s)-Mz(s))\right) B\sigma\left(B^\top[Pz(s) - M^*(w(s)-Mz(s))]\right)
$$ 
are also linearly bounded. Hence, applying Barbalat's Lemma and using Definition~\ref{def-cb}, we have
\begin{equation}
\lim_{t\rightarrow +\infty} z(t) = 0,
\qquad\: \lim_{t\rightarrow +\infty} B^\top M^* w(s) = 0.
\end{equation}
Let $\omega(z_0, w_0)$ be the $\omega$-limit set of the initial condition $(z_0, w_0)\in D(\cA)$, that is, is the set of all $(z^\star, w^\star)\in D(\cA)$ such that there exists an increasing sequence of time $(t_n)_{n\geq0}$ such that $z(t_n)\to z^\star$ and $w(t_n)\to w^\star$ in $H$ as $n$ goes to infinity.
Since $\cA$ is a $m$-dissipative operator, the trajectory $(z, w)$ is bounded in $D(\cA)$.
The injection of $D(\cA)$ in $X$ being compact, if $\omega(z_0, w_0)=\{(0, 0)\}$, then $(z, w)$ converges to the origin in the $X$-topology.
Since $(z, w)$ is bounded in time and since the injection of $D(\cA)$ in $X$ is compact, the positive orbit $\lbrace (z(t),w(t))\mid t\geq 0\rbrace$ is precompact in $X$. Therefore, according to the LaSalle's Invariance Principle for infinite-dimensional systems \cite[Theorem 3.1.]{slemrod1989mcss}, $\omega(z_0,w_0)$ is a non-empty compact set of $X$ that is invariant to the flow of \eqref{cl-system}. Hence, for any initial condition $(z_0^\star,w_0^\star)\in \omega(z_0,w_0)$, the corresponding solution $(z^\star,w^\star)\in C^1((0,\infty);X)\times C^0((0,\infty);D(\cA))$ of \eqref{cl-system} satisfies
\begin{equation}
\left\{
\begin{aligned}
&\frac{\d}{\dt} w^\star = S w^\star\\
&w^\star(0)=w_0^\star\\
&B^\top M^* w^\star = 0,\: z^\star = 0
\end{aligned}
\right.
\end{equation}
Since the pair $(S,B^\top M^*)$ is approximately observable in infinite-time, this implies that $w^\star_0=0$. Therefore, $\omega(z_0,w_0)=\lbrace (0,0)\rbrace$, which means that $(z,w)$ converges to $(0,0)$ in the $X$-topology. Moreover, according to Theorem~\ref{thm-wp}, $(0,0)$ is Lyapunov stable in the $X$-topology. Thus, $(0,0)$ is globally asymptotically stable, which concludes the proof of Theorem
~\ref{thm-GAS}.\qed
\end{proof}

\section{Observability Result for Skew-Adjoint Operators}
\label{sec:observability}

As mentioned earlier in Remark \ref{rem-pair},
we state now some sufficient conditions to verify the observability requirement of item \ref{pair-observable} in Theorem~\ref{thm-GAS}. We focus on the case in which $S$ is a skew-adjoint operator.
 
\color{black}
\begin{proposition}
\label{thm:observability}
Suppose that, in addition to Assumption~\ref{assumption1}, the following properties hold:
\begin{enumerate}[label=(\roman*),leftmargin=*, resume=hyp]
\item The operator $S$ is skew-adjoint and has compact resolvent;\label{item7}
\item The pair $(S^*,\Gamma^*)$ is approximately observable in infinite-time;\label{item8}
\item\label{item9} For all eigenvalue $\mu$ of $S$,
\begin{equation}
\label{non-resonance}
\rank
\begin{pmatrix}
A-\mu I_n & B\\
C & 0
\end{pmatrix}
=
n+p.
\end{equation}
\end{enumerate}
Then, the pair $(S, B^\top M^*)$ is approximately observable in infinite-time.
\end{proposition}

 \begin{remark}[About item~\ref{item3}.]
 Note that with item~\ref{item2}, i.e. $A$  Hurwitz, and \ref{item7}, i.e. $S$ skew-adjoint, item~\ref{item3}, i.e. the spectra of $A$ is $S$ are disjoint and nonempty, is trivially satisfied.
\end{remark}

\begin{proof}[Proposition \ref{thm:observability}]
Pick $\varphi\in D(S)$ an eigenvector of $S$ and denote $-\mu$ its corresponding eigenvalue. In order to apply \cite[Proposition 6.9.1]{tucsnak2009observation}, let us assume that $B^\top M^* \varphi = 0$ and look for a contradiction. Suppose that $\varphi\neq 0$.

Since $S$ is skew-adjoint, $\varphi$ is also an eigenvector of $S^*$, with corresponding eigenvalue $\mu$.
Then, according to \eqref{eq-Sylvester},
\begin{equation*}
\varphi^* SM - \varphi^* MA = -\varphi^* \Gamma C,
\qquad i.e.,
\qquad
\varphi^*M(\mu I_n-A)+\varphi^*\Gamma C = 0.
\end{equation*}
Hence,

\begin{equation}
\begin{pmatrix}
\varphi^*M
& \varphi^*\Gamma
\end{pmatrix}
\begin{pmatrix}
\mu I_n-A & B\\
C & 0
\end{pmatrix}
= 0.
\end{equation}

Since the matrix
$
\begin{pmatrix}
\mu I_n-A & B\\
C & 0
\end{pmatrix}
$
has full row rank, $\Gamma^*\varphi = 0$. According to \cite[Proposition 6.9.1]{tucsnak2009observation}, this contradicts the infinite-time observability of the pair $(S^*, \Gamma^*)$. This concludes the proof of Proposition \ref{thm:observability}.\qed 
\end{proof}

\begin{example}[The property \ref{pair-observable} for \eqref{hyperbolic}]
Consider again system \eqref{hyperbolic} introduced in Example~2. By some standard Sobolev injection theorems, it is easy to prove that the canonical embedding from $D(S)$ to $H$ is compact, i.e., $S$ has compact resolvent. Therefore, item \ref{item-resolvant} is satisfied. Moreover, assume that $M=\frac{N}{2}$ and $D_0$ and $D_1$ are orthogonal, which ensures that the corresponding operator $S$ is skew-adjoint, as noticed in \cite{russell1978controllability}. In \cite[Theorem 3.2]{russell1978controllability}, the exact controllability\footnote{See \cite[Definition 11.1.1.]{tucsnak2009observation} for a definition.} of the pair $(S,\Gamma)$ is proven, which, in turns, implies approximate observability in infinite-time of the pair $(S^*,\Gamma^*)$. Therefore, items~\ref{item7} and \ref{item8}
of Proposition~\ref{thm:observability} are satisfied. Hence, if $A$, $B$ and $C$ are given so that the conditions \ref{item9} is also satisfied, Proposition~\ref{thm:observability} does apply and hence also the condition~\ref{pair-observable} is satisfied, namely Theorem~\ref{thm-GAS} does apply.
\hfill$\triangle$
\end{example}

\color{black}

\section{Conclusions}\label{sec:conclusion}

In this paper, we have provided a methodology to globally stabilize a cascade system composed by a dissipative infinite-dimensional system and an ODE. This methodology relies on the so-called forwarding approach, which needs to introduce an infinite dimensional the Sylvester equation. Under an appropriate observability assumption, we prove that our design leads to a stabilizing feedback law. We also have provided a necessary condition implying this observability property. The proposed methodology can apply to other topics than the stabilization such as output regulation \cite{paunonen2010internal,paunonen2019stability} or the repetitive control \cite{weiss1999repetitive,astolfi2021repetitive}.

As further research lines, we may consider more complicated systems, such as \emph{nonlinear} infinite-dimensional dissipative systems. Let us mention for example the Korteweg-de Vries equation \cite{cerpa2013control}. In this case, we would not have to face with a Sylvester equation, which is a static equation, but with another dynamical nonlinear PDEs, for which the proof of the existence and the uniqueness of solution is not an easy task. Finally, another possible extension, can be to generalize the proposed forwarding approaches to different cascades of systems, for instance, in which each subsystem is  infinite-dimensional.


%
%

\bibliographystyle{spmpsci}      
\bibliography{bibsm}   

\end{document}